\documentclass{article}
\usepackage{amsmath, amsthm, amssymb}

\newtheoremstyle{theorem}
  {10pt}          
  {10pt}  
  {\sl}  
  {\parindent}     
  {\bf}  
  {. }    
  { }    
  {}     
\theoremstyle{theorem}
\newtheorem{theorem}{Theorem}

\newtheorem{lemma}[theorem]{Lemma}

\newtheoremstyle{defi}
  {10pt}          
  {10pt}  
  {\rm}  
  {\parindent}     
  {\bf}  
  {. }    
  { }    
  {}     
\theoremstyle{defi}
\newtheorem{definition}[theorem]{Definition}

\newtheorem{remark}[theorem]{Remark}

\begin{document}

\title{Constants of motion for fractional action-like
variational problems\footnote{Research Report CM06/I-24,
Department of Mathematics, University of Aveiro.}}

\author{Gast\~{a}o S. F. Frederico\\
\texttt{gfrederico@mat.ua.pt} \and
Delfim F. M. Torres\\
\texttt{delfim@mat.ua.pt}}
\date{Department of Mathematics\\
University of Aveiro\\
3810-193 Aveiro, Portugal}

\maketitle

\pagestyle{empty}
\thispagestyle{empty}

\begin{abstract}
We extend Noether's symmetry theorem to the fractional
Riemann-Liouville integral functionals
of the calculus of variations recently introduced by El-Nabulsi.

\bigskip


\smallskip

\end{abstract}


\section{Introduction}

The concept of \emph{symmetry} plays an important
role both in Physics and Mathematics.
Symmetries are described by transformations of the system,
which result in the same object after the transformation
is carried out. They are described mathematically
by parameter groups of transformations.
Their importance ranges from fundamental and theoretical aspects
to concrete applications, having profound implications
in the dynamical behavior of the systems, and
in their basic qualitative properties.

Another fundamental notion in Physics and Mathematics
is the one of \emph{constant of motion}.
Typical application of the constants of motion
in the calculus of variations
is to reduce the number of degrees of freedom,
thus reducing the problems to a lower dimension,
facilitating the integration of the differential
equations given by the necessary optimality conditions.

Emmy Noether was the first who proved, in 1918, that the
notions of symmetry and constant of motion are connected:
when a system exhibits a symmetry, then a constant of motion can be obtained.
One of the most important and well known illustrations of this
deep and rich relation, is given by the conservation of energy in Mechanics:
the autonomous Lagrangian $L(q,\dot{q})$,
correspondent to a mechanical system of conservative points,
is invariant under time-translations
(time-homogeneity symmetry), and
\begin{equation}
\label{eq:consEneg}
L\left(q,\dot{q}\right)
-\frac{\partial L}{\partial \dot{q}}\left(q,\dot{q}\right) \cdot \dot{q}
\equiv \text{constant}
\end{equation}
follows from Noether's theorem,
\textrm{i.e.}, the total energy of a conservative
system always remain constant in time, ``it cannot be created
or destroyed, but only transferred from one form into another''.
Expression \eqref{eq:consEneg}
is valid along all the Euler-Lagrange extremals $q$
of an autonomous problem of the calculus of variations.
The constant of motion \eqref{eq:consEneg} is known
in the calculus of variations as the 2nd Erdmann necessary condition;
in concrete applications, it gains different interpretations:
conservation of energy in Mechanics;
income-wealth law in Economics;
first law of Thermodynamics; etc.
The literature on Noether's theorem is vast,
and many extensions of the classical results of Emmy Noether
are now available in the literature (see \textrm{e.g.}
\cite{torresPortMath,Torres:CPAA:2004} and references therein).
Here we remark that constants of motion appear naturally in closed systems.

It turns out that in practical terms closed systems do not exist:
forces that do not store energy, so-called nonconservative
or dissipative forces, are always present in real systems.
Friction is an example of a nonconservative force.
Any friction-type force, like air resistance,
is a nonconservative force. Nonconservative forces remove energy from the systems
and, as a consequence, the constant of motion \eqref{eq:consEneg} is broken.
This explains, for instance, why the innumerable
``perpetual motion machines'' that have been proposed fail.
In presence of external nonconservative forces,
Noether's theorem and respective constants of motion cease to be valid.
However, it is still possible to obtain a Noether-type theorem
which covers both conservative (closed system)
and nonconservative cases \cite{CD:Djukic:1980,comGastaoParis06}.
Roughly speaking, one can prove
that Noether's conservation laws are still valid if a new term,
involving the nonconservative forces, is added to the standard
constants of motion.

The study of fractional problems of the calculus of variations
and respective Euler-Lagrange type equations
is a subject of strong current research because of its numerous applications:
see \textrm{e.g.} \cite{CD:Agrawal:2002,CD:BalAv:2004,CD:El-Na:2005,El-Nabulsi2005a,%
Klimek2001,MR1966935,Klimek2005,CD:Riewe:1996,CD:Riewe:1997}.
F.~Riewe \cite{CD:Riewe:1996,CD:Riewe:1997} obtained a version of
the Euler-Lagrange equations for problems of the calculus of
variations with fractional derivatives, that combines the
conservative and non-conservative cases. In 2002 O.~Agrawal
proved a formulation for variational problems with right and left
fractional derivatives in the Riemann-Liouville sense
\cite{CD:Agrawal:2002}. Then these Euler-Lagrange equations were
used by D.~Baleanu and T.~Avkar to investigate problems with
Lagrangians which are linear on the velocities
\cite{CD:BalAv:2004}. In \cite{Klimek2001,MR1966935} fractional problems of the calculus of
variations with symmetric fractional derivatives are considered and correspondent
Euler-Lagrange equations obtained, using both Lagrangian and Hamiltonian
formalisms. In all the above mentioned studies,
Euler-Lagrange equations depend on left and right
fractional derivatives, even when the problem depend only on
one type of them. In \cite{Klimek2005} problems depending on
symmetric derivatives are considered for which
Euler-Lagrange equations include only the
derivatives that appear in the formulation of the problem.
In \cite{CD:El-Na:2005,El-Nabulsi2005a} Riemann-Liouville
fractional integral functionals, depending on a parameter $\alpha$ but not on
fractional-order derivatives of order $\alpha$, are introduced and respective
fractional Euler-Lagrange type equations obtained.

A Noether-type theorem for problems of the calculus of variations
with fractional-order derivatives of order $\alpha$ is given in
\cite{FDA06}. Here we use the results of El-Nabulsi
\cite{CD:El-Na:2005,El-Nabulsi2005a} to prove a nonconservative
Noether's theorem in the new fractional action-like framework.


\section{Fractional action-like Noether's theorem}

We consider the fundamental problem of the calculus of variations
with Riemann-Liouville fractional integral, as considered by
El-Nabulsi \cite{CD:El-Na:2005,El-Nabulsi2005a}:
\begin{equation}
\label{Pi} I[q(\cdot)] = \frac{1}{\Gamma(\alpha)}\int_a^b
L\left(\theta,q(\theta),\dot{q}(\theta)\right)(t-\theta)^{\alpha-1}
d\theta \longrightarrow \min \, ,
\end{equation}
under given boundary conditions $q(a)=q_{a}$ and $q(b)=q_{b}$, where
$\dot{q} = \frac{dq}{d\theta}$, $\Gamma$ is the Euler gamma function,
$0<\alpha\leq 1$, $\theta$ is the intrinsic time, $t$ is the observer time,
$t \ne \theta$, and the Lagrangian $L :[a,b] \times \mathbb{R}^{n} \times
\mathbb{R}^{n} \rightarrow \mathbb{R}$ is a
$C^{2}$ function with respect to its arguments.
We will denote by $\partial_{i}L$ the partial derivative
of $L$ with respect to the $i$-th argument, $i = 1,2,3$.
Admissible functions $q(\cdot)$ are assumed to be $C^2$.

\begin{theorem}[\textrm{cf.} \cite{CD:El-Na:2005}]
\label{Thm:NonDtELeq} if $q$ is a minimizer of problem
\eqref{Pi}, then $q$ satisfies the following
\emph{Euler-Lagrange equation}:
\begin{equation}
\label{eq:elif}
\partial_{2} L\left(\theta,q(\theta),\dot{q}(\theta)\right)-\frac{d}{d\theta}
\partial_{3} L\left(\theta,q(\theta),\dot{q}(\theta)\right)
= \frac{1-\alpha}{t-\theta}\partial_{3}
L\left(\theta,q(\theta),\dot{q}(\theta)\right)\, .
\end{equation}
\end{theorem}

We now introduce the following definition of variational quasi-invariance up to a gauge term
(\textrm{cf.} \cite{torresPortMath}).

\begin{definition}[quasi-invariance of \eqref{Pi} up to a gauge term $\Lambda$]
\label{def:invaif} Functional \eqref{Pi} is said to be quasi-invariant under
the infinitesimal $\varepsilon$-parameter transformations
\begin{equation}
\label{eq:tinfif}
\begin{cases}
\bar{\theta} = \theta + \varepsilon\tau(\theta,q) + o(\varepsilon) \\
\bar{q}(\bar{\theta}) = q(\theta) + \varepsilon\xi(\theta,q) + o(\varepsilon) \\
\end{cases}
\end{equation}
up to the gauge term $\Lambda$ if, and only if,
\begin{multline}
\label{eq:invif1}
L\left(\bar{\theta},\bar{q}(\bar{\theta}),{\bar{q}}'(\bar{\theta})\right)
(t-\bar{\theta})^{\alpha-1}\frac{d\bar{\theta}}{d \theta} \\
= L\left(\theta,q(\theta),\dot{q}(\theta)\right)(t-\theta)^{\alpha-1}
+ \varepsilon (t-\theta)^{\alpha-1}\frac{d\Lambda}{d\theta}\left(\theta,q(\theta),
\dot{q}(\theta)\right) + o(\varepsilon)\, .
\end{multline}
\end{definition}

\begin{lemma}[necessary and sufficient condition for quasi-invariance]
\label{theo:cnsiif}
If functional \eqref{Pi} is quasi-invariant up to $\Lambda$ under
the infinitesimal transformations \eqref{eq:tinfif}, then
\begin{multline}
\label{eq:cnsiif}
\partial_{1}
L\left(\theta,q,\dot{q}\right)\tau+\partial_{2}
L\left(\theta,q,\dot{q}\right)\cdot\xi
+\partial_{3}
L\left(\theta,q,\dot{q}\right)\cdot\left(\dot{\xi}-\dot{q}\dot{\tau}\right)\\
+ L\left(\theta,q,\dot{q}\right)
\left( \dot{\tau} + \frac{1-\alpha}{t-\theta} \tau \right)
=\dot{\Lambda}\left(\theta,q,\dot{q}\right) \, .
\end{multline}
\end{lemma}

\begin{proof}
Equality \eqref{eq:invif1} is equivalent to
\begin{multline}
\label{eq:invif2}
\left[
L\left(\theta+\varepsilon\tau+o(\varepsilon),q+\varepsilon\xi+o(\varepsilon),
\frac{\dot{q}+ \varepsilon\dot{\xi}+o(\varepsilon)}{1+\varepsilon \dot{\tau}+o(\varepsilon)}\right)\right]
\left(t-\theta-\varepsilon\tau - o(\varepsilon)\right)^{\alpha-1}\left(1+\varepsilon \dot{\tau}+o(\varepsilon)\right)\\
= L\left(\theta,q,\dot{q}\right)(t-\theta)^{\alpha-1}
+ \varepsilon(t-\theta)^{\alpha-1}\frac{d}{d\theta} \Lambda\left(\theta,q,\dot{q}\right) +o(\varepsilon)
\, .
\end{multline}
Equation \eqref{eq:cnsiif} is obtained differentiating both sides
of equality \eqref{eq:invif2} with respect to
$\varepsilon$ and then putting $\varepsilon=0$.
\end{proof}

\begin{definition}[constant of motion]
\label{def:leico} A quantity $C\left(\theta,q(\theta),\dot{q}(\theta)\right)$,
$\theta \in [a,b]$, is said to be a \emph{constant of motion} if, and only if,
$\frac{d}{d\theta}C\left(\theta,q(\theta),\dot{q}(\theta)\right)=0$
for all the solutions $q$ of the Euler-Lagrange equation \eqref{eq:elif}.
\end{definition}

\begin{theorem}[Noether's theorem]
\label{theo:tnif} If the fractional integral \eqref{Pi} is quasi-invariant up to $\Lambda$,
in the sense of Definition~\ref{def:invaif}, and functions
$\tau(\theta,q)$ and $\xi(\theta,q)$ satisfy the condition
\begin{equation}
\label{eq:condif}
L\left(\theta,q,\dot{q}\right) \tau =
- \partial_{3} L\left(\theta,q,\dot{q}\right)\cdot(\xi - \dot{q}\tau) \, ,
\end{equation}
then
\begin{equation}
\label{eq:TeNet}
\partial_{3} L\left(\theta,q,\dot{q}\right)\cdot\xi(\theta,q)
+ \left[ L(\theta,q,\dot{q}) - \partial_{3}
L\left(\theta,q,\dot{q}\right) \cdot \dot{q} \right]
\tau(\theta,q)-\Lambda\left(\theta,q,\dot{q}\right)
\end{equation}
is a constant of motion.
\end{theorem}

\begin{remark}
Under our hypothesis \eqref{eq:condif} the necessary and sufficient
condition of quasi-invariance \eqref{eq:cnsiif} is reduced to
\begin{multline}
\label{eq:cnsiif1}
\partial_{1}
L\left(\theta,q,\dot{q}\right)\tau+\partial_{2}
L\left(\theta,q,\dot{q}\right)\cdot\xi +\partial_{3}
L\left(\theta,q,\dot{q}\right)\cdot\left(\dot{\xi}-\dot{q}\dot{\tau}\right) \\
+L\left(\theta,q,\dot{q}\right)\dot{\tau}-\frac{1-\alpha}{t-\theta}
\partial_{3}
L\left(\theta,q,\dot{q}\right)\cdot(\xi-\dot{q}\tau)=\dot{\Lambda}\left(\theta,q,\dot{q}\right)
\, .
\end{multline}
Conditions \eqref{eq:condif} and \eqref{eq:cnsiif1}
correspond to the generalized equations of Noether-Bessel-Hagen
of a non-conservative mechanical system \cite{CD:Djukic:1980}.
\end{remark}

\begin{proof}
We can write \eqref{eq:cnsiif1} in the form
\begin{multline}
\label{eq:cnsiif2}
\left[\partial_{1} L\left(\theta,q,\dot{q}\right)
+\frac{1-\alpha}{t-\theta}
\partial_{3} L\left(\theta,q,\dot{q}\right)\cdot \dot{q} \right] \tau
+ \left[ L\left(\theta,q,\dot{q}\right)
- \partial_3 L\left(\theta,q,\dot{q}\right)\cdot \dot{q} \right] \dot{\tau} \\
+ \left[ \partial_{2}
L\left(\theta,q,\dot{q}\right) - \frac{1-\alpha}{t
- \theta} \partial_3 L\left(\theta,q,\dot{q}\right)\right] \cdot \xi
+ \partial_3 L\left(\theta,q,\dot{q}\right) \cdot \dot{\xi} - \dot{\Lambda} = 0 \, .
\end{multline}
Using the Euler-Lagrange equation \eqref{eq:elif} equality \eqref{eq:cnsiif2}
is equivalent to
\begin{multline*}
\frac{d}{d\theta} \left[ L\left(\theta,q,\dot{q}\right)
- \partial_3 L\left(\theta,q,\dot{q}\right)\cdot \dot{q} \right] \tau
+ \left[ L\left(\theta,q,\dot{q}\right)
- \partial_3 L\left(\theta,q,\dot{q}\right)\cdot \dot{q} \right] \dot{\tau} \\
+ \frac{d}{d\theta} \left[ \partial_3 L\left(\theta,q,\dot{q}\right) \right] \cdot \xi
+ \partial_3 L\left(\theta,q,\dot{q}\right) \cdot \dot{\xi} - \dot{\Lambda} = 0
\end{multline*}
and the intended conclusion follows:
\begin{equation*}
\frac{d}{d\theta}\left[\partial_{3}
L\left(\theta,q,\dot{q}\right)\cdot\xi + \left(
L(\theta,q,\dot{q}) -
\partial_{3} L\left(\theta,q,\dot{q}\right) \cdot \dot{q} \right)
\tau-\Lambda\left(\theta,q,\dot{q}\right)\right]=0\, .
\end{equation*}
\end{proof}


\section{Examples}

In \cite[\S 4]{El-Nabulsi2005a} El-Nabulsi remarks
that conservation of momentum when $L$ is not a function
of $q$ or conservation of energy when $L$ has no explicit dependence
on time $\theta$ are no more true for a fractional order of integration
$\alpha$, $\alpha \ne 1$. As we shall see now, these facts are a trivial consequence
of our Theorem~\ref{theo:tnif}. Moreover, our Noether's theorem gives new explicit formulas
for the fractional constants of motion. For the particular case $\alpha = 1$ we recover
the classical constants of motion of momentum and energy.

Let us first consider an arbitrary fractional action-like problem
\eqref{Pi} with an autonomous $L$: $L\left(\theta,q,\dot{q}\right)
= L\left(q,\dot{q}\right)$. In this case $\partial_1 L = 0$,
and it is a simple exercise to check that \eqref{eq:cnsiif1}
is satisfied with $\tau = 1$, $\xi = 0$ and $\Lambda$ given by
\begin{equation*}
\dot{\Lambda} = \frac{1 - \alpha}{t - \theta} \,
\frac{\partial L}{\partial \dot{q}} \cdot \dot{q} \, .
\end{equation*}
It follows from our Noether's theorem (Theorem~\ref{theo:tnif}) that
\begin{equation}
\label{eq:ConsHam:alpha}
L\left(q,\dot{q}\right)
-\frac{\partial L}{\partial \dot{q}}\left(q,\dot{q}\right) \cdot \dot{q}
- (1 - \alpha) \int \frac{1}{t - \theta}
\frac{\partial L}{\partial \dot{q}}\left(q,\dot{q}\right) \cdot \dot{q} \, d\theta
\equiv \text{constant} \, .
\end{equation}
In the classical framework $\alpha = 1$ and
we then get from our expression \eqref{eq:ConsHam:alpha}
the well known constant of motion \eqref{eq:consEneg}, which
corresponds in mechanics to conservation of energy.

When $L$ is not a function of $q$ one has $\frac{\partial L}{\partial q} = 0$
and \eqref{eq:cnsiif1} holds true with $\tau = 0$, $\xi = 1$ and $\Lambda$ given by
\begin{equation*}
\dot{\Lambda} = - \frac{1 - \alpha}{t - \theta} \,
\frac{\partial L}{\partial \dot{q}}\left(\theta,\dot{q}\right)  \, .
\end{equation*}
The constant of motion \eqref{eq:TeNet} takes the form
\begin{equation}
\label{eq:consMomt:Frac}
\frac{\partial L}{\partial \dot{q}}\left(\theta,\dot{q}\right)
+ (1 - \alpha) \int \frac{1}{t - \theta}
\frac{\partial L}{\partial \dot{q}}\left(\theta,\dot{q}\right) d\theta \, .
\end{equation}
For $\alpha = 1$ \eqref{eq:consMomt:Frac} implies conservation of momentum:
$\frac{\partial L}{\partial \dot{q}} = const$.


\section*{Acknowledgments}

The first author acknowledges the support of
the \emph{Portuguese Institute for Development} (IPAD);
the second author the support by the \emph{Centre for Research on Optimization
and Control} (CEOC) from the ``Fundaç\~{a}o para a Ci\^{e}ncia e a Tecnologia'' (FCT),
cofinanced by the European Community Fund FEDER/POCTI.


\end{document}